\documentclass[12pt,letterpaper]{article}
\usepackage[utf8]{inputenc}
\usepackage[T1]{fontenc}
\usepackage[english]{babel}
\usepackage[margin=0.7in]{geometry}
\usepackage{microtype}
\usepackage{mathtools,amsfonts,amssymb}
\usepackage{appendix}
\usepackage{xcolor}
\allowdisplaybreaks
\usepackage[pdftex, colorlinks,  bookmarks=true, bookmarksnumbered=true, linkcolor=blue,  citecolor=red]{hyperref} 
\usepackage{natbib}
\usepackage{cleveref}
\usepackage[spanish, ruled,vlined]{algorithm2e}
\usepackage{amsthm}
\usepackage{subfigure} 
\usepackage{enumerate}
\usepackage{amsthm}     
\usepackage{multirow}

\def\diag{\textrm{diag }}

 \newtheorem{lemma}{Lemma}
\newtheorem{hypothesis}{Hypothesis}

\title{A variational model for wrapped phase denoising}
\author{Ivan May-Cen, Ricardo Legarda-Saenz, and Carlos Brito-Loeza}
\author{Ivan May-Cen, Ricardo Legarda-Saenz, and Carlos Brito-Loeza\thanks{\ All authors are with Computational Learning and Imaging Research group at Facultad de Matem\'aticas,\ Universidad Aut\'onoma de Yucat\'an,
		Mexico. \ Emails: ivan.mc@progreso.tecnm.mx; rlegarda, carlos.brito@correo.uady.mx.}  }
	\date{}
\begin{document}

\maketitle

\begin{abstract}
In this paper, we introduce a total variation based variational model for denoising wrapped phase images. Our model improves on former methods by preserving discontinuities of the phase map and enforcing the fundamental Pythagorean trigonometric identity between the real and imaginary parts of the phase map enhancing the quality of the restored phase. The existence and uniqueness of the solution of our model is proven using standard methods. Further, we provide a fast fixed point method for finding the numerical solution and prove its  convergence. Experiments on both synthetic and real patterns verify our findings.

\end{abstract}

\section{Introduction}
Optical measurement technology is a technique for estimating through light object's properties such as shape and deformation in a contact-free process. This technology widely employed for quality assurance in the manufacturing industry demands reliable and accurate methods. By extracting the local wrapped phase from one or a collection of interference fringe pattern images and then unwrapping them, one can estimate some physical properties of the objects being measured \citep{kulkarni2017single, servin2014fringe, surrel2000fringe}.  In interferometric synthetic aperture radar (InSAR) the acquired wrapped phase is the principal value of the absolute phase, i.e., a modulo-$2\pi$ observation.

Either, optical metrology or InSAR, the presence of noise  in the wrapped phase data,  mainly due to capture devices \citep{warlick2015errors, wyant2003dynamic, huntley1997random}, is the main cause of inaccuracies. Due to this, a denoising pre-processing step is usually performed on the wrapped phase map before unwrapping it. In this work, we introduce a variational model specifically designed for denoising fringe pattern  images that preserves the local wrapped phase content. We argue that general denoising algorithms already in use remove relevant information from the wrapped phase images, decimating the quality of the final measurements.

The mathematical representation of a phase map is 
\begin{equation}\label{ec-exp-fase}
U = A \exp (i \phi)
\end{equation}
where $A$ is the magnitude and $\phi$ the phase map. Most methods assume that the wrapped phase (modulo-$2\pi$ ) is given by
\begin{equation}\label{ec-atan2}
\psi = \arctan2(\sin \phi, \cos \phi).
\end{equation}
where the function $\arctan2$ is defined as
\begin{gather}
	 \arctan2(y,x) = \left\{ \begin{array}{ll}
		\arctan(y/x) & \mbox{if $x > 0$};\\
		\arctan(y/x) +  \pi & \mbox{if $x < 0$ and $y \geq 0$};\\
		\arctan(y/x) - \pi & \mbox{if $x < 0$ and $y < 0$};\\
		+\pi/2 & \mbox{if $x = 0$ and $y > 0$};\\
		-\pi/2 & \mbox{if $x = 0$ and $y < 0$};\\		
		\text{undefined} & \mbox{if $x = 0$ and $y = 0$}.\end{array} \right.
\end{gather}

The above assumption, implies the need to preserve discontinuities induced by the $\arctan$ function and discontinuities that $\phi$ may already have. Therefore, denoising algorithms for wrapped phase maps must to be capable of preserving both types of discontinuities. Up to our knowledge, no method has been proposed that considers this fact. 

Further, denoising independently the noisy real and imaginary (sine and cosine) images, as suggested in \cite{strobel1996processing}, will lead to inaccuracies when reconstructing the wrapped phase by means of the arctan2 function. To avoid this, we introduce a term in our model that couples both trigonometric functions improving on the reconstruction of $\psi$.

In what follows, a review of some methods that have been used for phase denoising is presented hereunder.

In  \cite{kulkarni2018phase, kulkarni2020fringe}, the authors proposed a simultaneous phase unwrapping and noise filtering algorithms with a probabilistic approach based on random Markov fields and Kalman filters; their proposals work well for continuous phases, but nothing is mentioned about discontinuous phase maps. \cite{li2017general} proposed a phase filtering method based on a variational decomposition model with total variation regularization; they did not show any results with discontinuous phase maps either. \cite{medina2017filtering} presented a global filtering process that makes use of the local frequencies of the wrapped phase map and tunes each pixel to its instantaneous frequency; this linear filtering method only works for continuous phases. \cite{kemao2008windowed, kemao2010windowed} proposed different methods for discontinuous maps, however they require assuming the phase to be locally polynomial and the overall process requires quality windows which adds complexity to the scheme and uncertainty if the size of the window is not chosen right. Modern methods for noise filtering based on deep learning algorithms  \citep{yan2020wrapped, zhang2019phase,yan2019fringe, yan2020deep} only perform well for the continuous phase maps.

In this paper we adopt the variational approach which in many other image processing techniques has proven to perform well by preserving discontinuities and noise filtering, see for instance \citep{sartor2008variational, vese2016variational, song2017insar} and references therein.  Here, a variational model for denoising wrapped phase maps that not only preserves discontinuities but also the Pythagorean relation between the real and imaginary parts of the phase map is proposed and analyzed. The structure of this work is as follows: in Section 2 we introduce our proposed variational model; Section 3 is dedicated to the analysis of the model and the derivation of the Euler-Lagrange equations; in Section 4, we introduce techniques for the numerical solution of our model. In particular, we present a fixed point method, its proof of convergence  and computational realization; in Section 5 we show some experiments to illustrate the performance of our model and compare with state of the art methods for denoising. Finally, in Section 6, we present our conclusions.

\section{A variational model for denoising wrapped phase maps}

In this section, we introduce our denoising model for wrapped phase images. 

The mathematical model of a wrapped phase map $\hat{\psi}$ polluted with additive Gaussian noise $\eta$ is well known and given by 
\begin{equation}
	\hat{\psi} = \psi + \eta
\end{equation}
where $\psi$ is the true wrapped phase map. 

Maybe the most used technique for filtering wrapped phase maps is the one proposed in \cite{strobel1996processing} where $\hat U_{real}= A\cos\hat\psi$ and $\hat U_{im}= A\sin\hat\psi$ are filtered separately by standard methods obtaining estimations $\bar{U}_{real}$ and $\bar{U}_{im}$. Then, the estimated $\bar{\psi}$ to $\psi$ is recovered by $\bar{\psi}=\arctan(\bar{U}_{im}/\bar{U}_{real})$ and $\bar A=[(\bar{U}_{real})^2+(\bar{U}_{im})^2]^{1/2}$. The drawback of this method is that modulo-$2\pi$ discontinuities are likely to be lost in the process due to the periodicity of sine and cosine functions.
		
Since there is no problem in computing the magnitude $A$ for an arbitrary phase map and normalize it, in this work, without loss of generality, we consider the case $A=1$ which yields the simplification $\hat U_{real}= \cos\hat\psi$ and $\hat U_{im}= \sin\hat\psi$.  Then, we apply the following minimization
\begin{gather*}
(\bar U_{real}, \bar U_{im}) =	\arg\min_{(x,y)} F(x,y),
\end{gather*}

where 
\begin{equation}\label{ec-funcional}
\begin{array}{rcl}
F & \equiv & \displaystyle{\frac{\lambda_1}{2}\int_{\Omega}(\bar U_{real}-\hat U_{real})^2d\Omega
+ \frac{\lambda_2}{2}\int_{\Omega}(\bar U_{im}-\hat U_{im})^2d\Omega
+ \frac{\lambda_3}{2}\int_{\Omega}((\bar U_{real})^2 + (\bar U_{im})^2-1)^2d\Omega}\\
& & + \displaystyle{\int_{\Omega}|\nabla \bar U_{real}|d\Omega
+ \int_{\Omega}|\nabla \bar U_{im}|d\Omega}
\end{array}
\end{equation}
where $\Omega\subset \mathbb{R}^2$ is the integration domain, $\lambda_1, \lambda_2$, and $\lambda_3$ are regularization parameters. The first two terms in $F$ constitute the fitting function following Str\"obel's technique for filtering wrapped phase maps \cite{strobel1996processing}. The third term is a main contribution of our model and forces the fundamental Pythagorean identity $(\bar U_{im})^2 + (\bar U_{real})^2=1$ to be satisfied improving on the estimation of $\psi$ by means of  the $\arctan$ function. Finally, the last two terms are there to regularize the solutions $\bar U_{im}, \bar U_{real}$. The total variation regularizer is well known to allow for discontinuities in the solution and to removing noise at the same time \citep{chan2005aspects}.

\section{Model analysis}

In this section, we prove existence and uniqueness of the solution of our model using standard techniques. 
Let us start by assuming that $\Omega$ is a convex region and defining the bounded variation norm in $\Omega$ by $||\cdot||_{BV}$.

\subsubsection*{Existence}
We start by showing that $F$ is BV-coercive. To this end, note that
\begin{gather}
	 \int_{\Omega}(\bar U_{im}-\hat U_{im})^2d\Omega \geq 0, \label{p1} \\\
		\int_{\Omega}(\bar U_{real}-\hat U_{real})^2d\Omega \geq 0, \label{p2} \\
	 \int_{\Omega}((\bar U_{real})^2 + (\bar U_{im})^2-1)^2d\Omega \geq 0. \label{p3}
\end{gather}
Then, it is immediate that 
\begin{gather}
	\label{coercive}
	\lim_{||(\bar U_{real}, \bar U_{im}) ||_{BV} \rightarrow \infty }  F(\bar U_{real},\bar U_{im}) = \infty
\end{gather}
meaning $F$ is BV-coercive. 

Weakly lower-semi-continuity of (\ref{p1}), (\ref{p2}) and, (\ref{p3}) comes from the weakly lower semi-continuity of the norms on Banach spaces. The same property holds for $||\bar U_{real}||_{BV}$ and $||\bar U_{im}||_{BV}$ as proven in \cite{acar1994analysis}. Finally, with $F$ being weakly lower semi-continuous and BV-coercive, the existence of the solution in the BV space is guaranteed following \cite{acar1994analysis}. 

\subsubsection*{Uniqueness}
Uniqueness of the solution comes from the model being convex. Convexity of (\ref{p1}),  (\ref{p2}), $||\bar U_{real}||_{BV}$, and $||\bar U_{im}||_{BV}$ are proven in \cite{acar1994analysis}.  Therefore just rests to prove convexity of (\ref{p3}). 
In Section \ref{fps} we show that (\ref{p3}) can be seen as $||Ku-d||_2^2$ with $u=(\bar U_{real}; \bar U_{im})$ and $K$ a bounded operator. It is immediate to see that (\ref{p3}) is therefore convex \cite{boyd2004convex}.

Note also that $F$ is bounded below by zero.
\subsection{Euler-Lagrange equations}
In this section, we derive the Euler-Lagrange equations of the model. Let us define $\varphi \in C(\Omega)$ with compact support and $\varepsilon$ a given scalar. Then, we set the first variations equal to zero, i.e
\begin{gather*}
\delta F_{real} \equiv \frac{d}{d\varepsilon}F(\bar U_{real}+\varepsilon\varphi, \bar U_{im})\Big|_{\varepsilon=0} =0, \;\; \text{and}\;\; \\
\delta F_{im} \equiv \frac{d}{d\varepsilon}F(\bar U_{real}, \bar U_{im}+\varepsilon\varphi)\Big|_{\varepsilon=0}=0.
\end{gather*}
We will detail the process for $\bar U_{real}$ only since for $\bar U_{im}$ is very similar
\begin{align}
\delta F_{real}& =  \frac{d}{d\varepsilon} \left\{ \frac{\lambda_1}{2}\int_{\Omega} (\bar U_{real}+\varepsilon\varphi - \hat U_{real})^2d\Omega 
 +\frac{\lambda_2}{2}\int_{\Omega} (\bar U_{im} - \hat U_{im})^2 d\Omega   \right.\\
& + \frac{\lambda_3}{2}\int_{\Omega} ((\bar U_{real}+\varepsilon\varphi)^2 +(\bar U_{im})^2 - 1)^2d\Omega\\
& \left.+ \int_{\Omega} |\nabla(\bar U_{real}+\varepsilon\varphi)|d\Omega + 
 	\int_{\Omega} |\nabla \bar U_{im}|d\Omega\right\}\Big|_{\varepsilon=0} \\
 &= \lambda_1 \int_{\Omega} (\bar U_{real}-\hat U_{real})\varphi \; d\Omega 
 + 2\lambda_3 \int_{\Omega}\left((\bar U_{real})^2+(\bar U_{im})^2-1 \right)\bar U_{real} \varphi \;d\Omega \\ &+ \int_{\Omega} \frac{\nabla \bar U_{real}}{|\nabla \bar U_{real}|}\cdot \nabla\varphi \; d\Omega.
\end{align}
Now, applying Green's Theorem in the last term,
\begin{equation}
\int_{\Omega} \frac{\nabla \bar U_{real}}{|\nabla \bar U_{real}|}\cdot \nabla\varphi d\Omega = -\int_{\Omega} \nabla \cdot \left( \frac{\nabla \bar U_{real}}{|\nabla \bar U_{real}|}\right)\varphi d\Omega +\oint_{\partial\Omega} \left(\frac{\nabla \bar U_{real}}{|\nabla \bar U_{real}|}\cdot \mathbf{n} \right)\varphi
d(\partial\Omega) 
\end{equation}
we get the following nonlinear partial differential equation (PDE)
\begin{equation}
	\label{pdels}
-\nabla \cdot \left( \frac{\nabla \bar U_{real}}{|\nabla \bar U_{real}|} \right) +  \lambda_1(\bar U_{real}-\hat U_{real})+2\lambda_3\left((\bar U_{real})^2+(\bar U_{im})^2-1\right)\bar U_{real} = 0  \;\;\;\text{in}\; \Omega
\end{equation}
with Neumann's boundary condition
\begin{gather}
	\frac{\nabla \bar U_{real}}{|\nabla \bar U_{real}|}\cdot \mathbf{n}=0 \;\;\text{on}\; \partial \Omega
\end{gather}
where $\mathbf{n}$ denotes the outward unit normal vector on $\partial \Omega$.

The correspondent PDE for $\bar U_{im}$ is
\begin{equation}
	\label{lcpde}
	-\nabla \cdot \left( \frac{\nabla \bar U_{im}}{|\nabla \bar U_{im}|} \right) +  \lambda_2(\bar U_{im}-\hat U_{im})+2\lambda_3\left((\bar U_{real})^2+(\bar U_{im})^2-1\right)\bar U_{im} = 0  \;\;\;\text{in}\; \Omega
\end{equation}
with Neumann's boundary condition
\begin{gather}
	\frac{\nabla \bar U_{im}}{|\nabla \bar U_{im}|}\cdot \mathbf{n}=0 \;\;\text{on}\; \partial \Omega.
\end{gather}

\section{Numerical solution}
The nonlinear PDEs in (\ref{pdels}) and (\ref{lcpde}) are second order nonlinear and anisotropic equations, therefore it is expected ordinary numerical methods to find it hard to converge quickly to the solution. This is the case of the  popular gradient descent method which due to stability restrictions maybe very slow to converge to the solution. In both PDEs, when computing the discrete diffusion coefficient $D=\frac{1}{|\nabla  \bar U_{real}|}$ for  (\ref{pdels}), it has to be regularized to avoid division by zero by adding a vary small scalar $\beta > 0$ i.e. $D_\beta=\frac{1}{|\nabla  \bar U_{real} + \beta|}$ . In \citep{marquina2000explicit}, it was showed that as $\beta \rightarrow 0$, the time-step parameter in the gradient descent scheme must tend to zero as well to guarantee stability of the numerical method. Of course, the smaller $\beta$ is, the larger number of iterations are required. 

Other numerical methods that have been developed for similar PDEs but still have to be tested in our wrapped phase denoising model are the dual formulations in \citep{chambolle2004algorithm, chan1999nonlinear}, alternating minimization methods \citep{huang2008fast} and references therein, Newton's based method  \citep{ng2007semismooth}, Bregman iterative algorithms \citep{osher2005iterative} and fixed point methods \citep{vogel1996iterative}. 
In this work, we present an efficient fixed point method with fast convergence and leave the rest for the future work.

\subsection{Fixed point method}
Here, we will introduce our fixed point algorithm for the solution of PDEs  (\ref{pdels}) and  (\ref{lcpde}). At first instance, the following fixed point scheme should work for $\bar U_{real}$ with $\bar U_{im}$ fixed
\begin{equation}
	\label{fppde1}
	\left(-\nabla \cdot \left( \frac{\nabla }{|\nabla \bar U_{real}^k|} \right) +  \lambda_1  +2\lambda_3\left((\bar U_{real}^k)^2+(\bar U_{im}^k)^2-1\right) \right) \bar U_{real} ^{k+1}=   \lambda_1  \hat U_{real}  \;\;\;\text{for}\; k=1,2,3,\ldots
\end{equation}

However, note that the term $(\bar U_{real})^2+(\bar U_{im})^2-1 $ may be negative at some iteration and for $2\lambda_3 \geq \lambda_1$ the resultant matrix is no longer strictly diagonally dominant. Therefore, algorithms such as Gauss-Seidel and Jacobi are not guaranteed to converge.
To fix this problem, we move one term to the right hand side resulting in
\begin{equation}
	\label{fp1}
	\left(-\nabla \cdot \left( \frac{\nabla }{|\nabla \bar U_{real}^k|} \right) +  \lambda_1 +  2\lambda_3\left((\bar U_{real}^k)^2+(\bar U_{im}^k)^2 \right)  \right) \bar U_{real}^{k+1} =   \lambda_1  \hat U_{real} + 2\lambda_3\bar U_{real}^k \;\;\;\text{in}\; \Omega.
\end{equation}
 Now, $ \lambda_1 +  2\lambda_3\left((\bar U_{real}^k)^2+(\bar U_{im}^k)^2 \right)$ which is always positive, is added to the main diagonal, and we have strict diagonal dominance guaranteed.
Following the same procedure, we can get a fixed point scheme for  $\bar U_{im}$ with $\bar U_{real}$ fixed
\begin{equation}
	\label{fp2}
	\left(-\nabla \cdot \left( \frac{\nabla }{|\nabla \bar U_{im}^k|} \right) +  \lambda_2 +  2\lambda_3\left((\bar U_{real}^k)^2+(\bar U_{im}^k)^2 \right)  \right) \bar U_{im}^{k+1} =   \lambda_2  \hat U_{im} + 2\lambda_3\bar U_{im}^k \;\;\;\text{in}\; \Omega.
\end{equation}
In short notation, we can represent (\ref{fp1}) and (\ref{fp2}) as
\begin{gather}
	H(\bar U_{real}^k) \bar U_{real}^{k+1}= f(\bar U_{real}^k), \label{fp1a}\\
	H(\bar U_{im}^k) \bar U_{im}^{k+1}= f(\bar U_{im}^k) \label{fp1b}
\end{gather}
and solve in alternating way first for $\bar U_{real}$ with $\bar U_{im}$ fixed and then the opposite way. However we can do better and create the system
\begin{gather}
	\label{fpv}
	L(u^{k}) u^{k+1} = c^{k}  
\end{gather}
where 
\begin{gather*}
	L(u^{k}) = \left( \begin{array}{cc}H(\bar U_{real}^{k})& 0\\
		0 & H(\bar U_{im}^{k}) \end{array}\right), 
	\\u^{k+1}=\left( \begin{array}{c}\bar U_{real}^{k+1}\\ \bar U_{im}^{k+1}\end{array}\right), 
	\\c^{k} = \left( \begin{array}{c} f(\bar U_{real}^k)\\  f(\bar U_{im}^k)\end{array}\right).
\end{gather*}

We note that the matrix $L(u)$ inherits the same nice properties of $H$ i.e., $L$ is symmetric, positive definite and strictly diagonal dominant.

\subsection{Fixed point algorithm convergence}\label{fps}

To prove convergence of our proposed fixed point scheme (\ref{fpv}), we will use  Weiszfeld method, see \citep{chan1999convergence, shi2008convergence,brito2020fast}. Here we use a lexicographic ordering on $ \bar U_{real} $ and $ \bar U_{im}$, both in $\mathbb{R}^{N}$, to define $u \in \mathbb{R}^{2N}$ as the column vector $u=( \bar U_{real}, \bar U_{im})$, where  $N= m \times n$. Then, our model in (\ref{ec-funcional}) can be written as follows
\begin{gather}
	\label{min2}
	\arg\min_u F(u)
\end{gather}
where
\begin{gather}\label{ec-F}
	F(u)  =   \sum_l |A_{l}^T u|_{\beta}  +\frac{\lambda_1}{2}||u - d_1||_2^2 
		+ \frac{\lambda_2}{2}||u - d_2||_2^2  +\frac{\lambda_3}{2}||Ku - d_3||_2^2,
\end{gather}
$A^T_l \in \mathbb{R}^{2\times 2N}$ is an array that contains the gradient of $u$ that is approximated using finite differences. Here, $l$ defined as 
\begin{gather}\label{ldef}
	l = \left\{ \begin{array}{ll}
		l=(i-1)n+j & \mbox{if $(i,j)$ belongs to $\bar U_{real}$}, \\
			l=N+(i-1)n+j & \mbox{if $(i,j)$ belongs to $\bar U_{im}$}\end{array} \right. 
\end{gather}
stands for the pixel coordinates where the gradient is computed and $|A^T_l u|_{\beta} = \sqrt{ (\partial_x u_l)^2 + (\partial_y u_l)^2 +\beta}$ where $\beta>0$ is a regularization parameter. 

Note that $A^T_l u$ can be split in two as $A^T_l u = (A^T_{1l}, 0)u + (0, A^T_{2l})u$. 
For instance, the operator $(A^T_{1l}, 0)$ computes the gradient of the first half of $u$, i.e. $\bar U_{real}$ and pads the vector with zeros. We will name $B_{1l}=(A^T_{1l}, 0)$ and $B_{2l}=(0, A^T_{2l})$.
Then, the two regularization terms in (\ref{ec-funcional})  can be compactly expressed as $ \sum_l |A_{l}^T u|_{\beta} $.

Further, the set $\{d_1, d_2, d_3\}$ has elements in $\mathbb{R}^{2N}$, $K \in \mathbb{R}^{2N\times 2N}$ is a bounded operator with $K^*$ its adjoint. 
\begin{gather}\label{ec-def-terminos}
	K = \left(\begin{array}{cc}
		\bar U_{real}^k & \bar U_{im}^k\\
		-\bar U_{im}^k & \bar U_{real}^k
	\end{array}
	\right) , 
	d_1 = \left( \begin{array}{c}
		\hat U_{real}\\
		0 
	\end{array}\right),
	d_2 = \left( \begin{array}{c}
		0\\
		\hat U_{im} 
	\end{array}\right), 
	d_3 = \left( \begin{array}{c}
	1\\
	0
\end{array}\right), 
\end{gather}
with $1$ the vector with all entries equal to one and $0$ the null vector.

Note that the fixed point algorithm in (\ref{fpv}) then can be written as
\begin{equation}\label{ec-forma}
	\sum_l A_l\Big( \frac{A^T_l u^{k+1}}{|A^T_l u^k|_{\beta}}\Big) +  \lambda_1 (u^{k+1}-d_1)+\lambda_2 (u^{k+1}-d_2)+2\lambda_3 K^*(Ku^{k+1}-d_3)=0.
\end{equation}

\subsubsection{Generalized Weiszfeld's method}
The generalized Weiszfeld’s method consists of choosing a uniformly strictly convex quadratic function $G(w,u)$ that approximates $F(u)$ under the assumptions of Hypothesis \ref{hipotesis} below

\begin{hypothesis}\label{hipotesis}
\hspace{0.1cm}
\begin{enumerate}
\item $G(w,u) = F(u) + (w-u, F'(u))+\frac{1}{2}(w-u, C(u)(w-u))$.
\item $C(u)$ is continuous.
\item $\lambda_{\min}(C(u)) \geq \mu >0\ \forall u$.
\item $F(w)\leq G(w,u)\ \forall w$.
\end{enumerate}
\end{hypothesis}

Weiszfeld generalized method solves (\ref{min2}) by the following iterative scheme
\begin{equation}
	u^{k+1} = \min_{w} G(w,u^k)
\end{equation}
where $G(w, u^{k})$ is a quadratic approximation of $F(u)$ at iteration $k$. 
Under Hypothesis \ref{hipotesis} and for fixed $u$, $G(w,u)$ is coercive, bounded and strictly convex. Therefore the minimum exists and can be computed by 
\begin{equation}\label{ec-vec-PF-G}
0 = G'_w(u^{k+1}, u^k) = F'(u^k)+C(u^k)(u^{k+1}-u^k)
\end{equation}
with global and linear convergence \cite{chan-mulet, shi2008convergence}.

In what follows, we will show with the help of the following Lemma, that (\ref{ec-forma}) and  (\ref{ec-vec-PF-G}) are equivalent. 

\begin{lemma}\label{lema-C}
If $C(u)$ is given by
\begin{equation}
C(u) = 
B_{1} \diag \left( \frac{1}{|B^T_{1l}|_\beta}I_2 \right) B^T_{1}
 + 
B_{2} \diag \left( \frac{1}{|B^T_{2l}|_\beta}I_2 \right) B^T_{2}
+ 
 \lambda_1 + \lambda_2+ 2\lambda_3K^*K
\end{equation}
for $B_1=(B_{11}, \ldots, B_{1m})$, $I_2$ the $2 \times 2$ identity matrix and $B_2$ defined in a similar way. Then 
\begin{equation}
G(w,u) = F(u)+(w-u, F'(u))+\frac{1}{2}(w-u, C(u)(w-u))
\end{equation}
defines a generalized Weiszfeld method and (\ref{ec-forma})  is equivalent to (\ref{ec-vec-PF-G}).
\end{lemma}
\begin{proof}
The proof is similar to that one given for Lemma 6.1 in \cite{shi2008convergence}. First, we show that  property 1 in Hypothesis \ref{hipotesis}  is satisfied, i.e that (\ref{ec-forma}) and (\ref{ec-vec-PF-G}) are equivalent.

\begin{align}
F'(u^k)+C(u^k)(u^{k+1}-u^k) &= \lambda_1 (u^k-d_1)+\lambda_2 (u^k-d_2)+2\lambda_3 K^*(Ku^k-d_3) \nonumber\\
  &+\displaystyle{\sum_l B_{1l}\Bigg( \frac{B_{1l}^Tu^k}{|B_{1l}^Tu^k|_{\beta}}  \Bigg)}
  +\displaystyle{\sum_l B_{2l}\Bigg( \frac{B_{2l}^Tu^k}{|B_{2l}^Tu^k|_{\beta}}  \Bigg)} \nonumber \\
 &+ \displaystyle{\sum_l B_{1l} \Bigg( \frac{B_{1l}^T u^{k+1}}{|B_{1l}^Tu^k|_{\beta}} \Bigg) }
  +\displaystyle{\sum_l B_{2l}\Bigg( \frac{B_{2l}^T u^{k+1}}{|B_{2l}^Tu^k|_{\beta}} \Bigg) } \nonumber \\
  &+\lambda_1 u^{k+1}+\lambda_2 u^{k+1}+2\lambda_3 K^*Ku^{k+1} \nonumber \\
  &-\displaystyle{\sum_l B_{1l}\Bigg( \frac{B_{1l}^T u^k}{|B_{1l}^Tu^k|_{\beta}} \Bigg) }
  -\displaystyle{\sum_l B_{2l}\Bigg( \frac{B_{2l}^T u^k}{|B_{2l}^Tu^k|_{\beta}} \Bigg) }\nonumber \\
  &-\lambda_1 u^{k}-\lambda_2 u^{k}-2\lambda_3 K^*Ku^{k} \nonumber\\
  &=  \lambda_1 (u^{k+1}-d_1)+\lambda_2 (u^{k+1}-d_2)+2\lambda_3 K^*(Ku^{k+1}-d_3) \nonumber \\
     &+ \displaystyle{\sum_l B_{1l}\Bigg( \frac{B_{1l}^Tu^{k+1}}{|B_{1l}^Tu^k|_{\beta}}  \Bigg)} \nonumber 
  +\displaystyle{\sum_l B_{2l}\Bigg( \frac{B_{2l}^Tu^{k+1}}{|B_{2l}^Tu^k|_{\beta}}  \Bigg)} \\
  &=  \lambda_1 (u^{k+1}-d_1)+\lambda_2 (u^{k+1}-d_2)+2\lambda_3 K^*(Ku^{k+1}-d_3)  \nonumber  \\
  &+\sum_l A_l\Big( \frac{A^T_l u^{k+1}}{|A^T_l u^k|_{\beta}}\Big).
\end{align}

Property 2 holds because of the continuity of $C(u)$ is guaranteed by $|B_{1l}^Tu^k|_\beta>0$ and $|B_{2l}^Tu^k|_\beta>0$.

Property 3 is fulfilled by noticing that $K$ is a bounded nonzero operator.

Property 4, is straightforward from the procedure in the proof of Lemma 6.1 in \cite{shi2008convergence} and will be omitted here. 

\end{proof}


\subsection{Numerical realization}
Here we show the numerical realization of the fixed point method for the Euler-Lagrange equations. To this end, we define the discrete domain as $\Omega_h=[0,m]\times[0,n]$ and without loss of generalization we set the spatial step size $h=h_x = h_y=1$ over the grid defined as  $\Gamma_h = \{(x,y) : x=ih_x, y=jh_y, i, j \in \mathbb{Z} \}$.

To approximate the partial derivatives of a given variable $w$, we use the finite differences method. 
The norm of the gradient of $w$ at points $({i+1,j} ), ({i-1,j}), ({i,j-1})$, and $({i,j-1})$ is computed as follows
\begin{gather}
|\nabla w|_{i+1,j}  =  \sqrt{\frac{(w_{i+1,j}-w_{i,j})^2}{h_x^2} + \frac{(w_{i,j+1}-w_{i,j})^2}{h_y^2} },\\
|\nabla w|_{i-1,j}  =  \sqrt{ \frac{(w_{i,j}-w_{i-1,j})^2}{h_x^2} +  \frac{(w_{i-1,j+1}-w_{i-1,j})^2}{h_y^2} },\\
|\nabla w|_{i,j+1}  =\sqrt{\frac{(w_{i+1,j}-w_{i,j})^2}{h_x^2} +  \frac{(w_{i,j+1}-w_{i,j})^2}{h_y^2}  },\\
|\nabla w|_{i,j-1}  =  \sqrt{ \frac{(w_{i+1,j-1}-w_{i,j-1})^2}{h_x^2} +  \frac{(w_{i,j}-w_{i,j-1})^2}{h_y^2} }.
\end{gather}
Then, the curvature term $\nabla \cdot \frac{\nabla u}{|\nabla u|} $ is approximated at point $({i,j})$ by 
\begin{gather}
\nabla \cdot \frac{\nabla w}{\sqrt{|\nabla w|^2+\beta}}  =  \frac{(w_{i+1,j}-w_{i,j})/h_x}{\sqrt{|\nabla w|_{i+1,j}^2 + \beta}} - \frac{(w_{i,j}-w_{i-1,j})/h_x}{\sqrt{|\nabla w|_{i-1,j}^2+\beta}} + \frac{(w_{i,j+1}-w_{i,j})/h_y}{\sqrt{|\nabla w|_{i,j+1}^2+\beta}} - \frac{(w_{i,j}-w_{i,j-1})/h_y}{\sqrt{|\nabla w|_{i,j-1}+\beta}}
\end{gather}
where $\beta>0$ is a regularization parameter to avoid division by zero.
Finally, the boundary conditions are enforced by 
\begin{equation}\label{ec-cond-frontera-neumann}
w_{i,0} = w_{i,1},\  w_{i,n+1} = w_{i,n},\  w_{0,j} = w_{1,j}, \  w_{m+1,j} = w_{m,j}.
\end{equation}

To realize the fixed point iteration, first define $\sqrt{|\nabla \bar U_{real}|_{i,j}^2 + \beta} = G_{i,j}$, and $\sqrt{|\nabla \bar U_{im}|_{i,j}+\beta} = \bar G_{i,j}$. Then, by taking advantage that the resulting matrix is strictly diagonal dominant, each iteration of the fixed point method can be solved using the Gauss Seidel method as follows:

\begin{gather}\label{ec-PF-Is-Ic}
\bar u^{k+1}_{l}	= \left( \begin{array}{c}
		(\bar U_{real}^{k+1})_l    \\
		(\bar U_{im}^{k+1} )_l \end{array}  \right)
	=
	\left( \begin{array}{cc}
	(\mathcal{N}_{real})_{i,j}/(\mathcal{D}_{real})_{i,j} \\
	(\mathcal{N}_{im})_{i,j}/(\mathcal{D}_{im})_{i,j}
	 \end{array}  \right)
\end{gather}
with $l$ as in (\ref{ldef}) and 
\begin{gather}
 (\mathcal{N}_{real})_{i,j}=  G^k_{i+1,j} (\bar U_{real}^k)_{i+1,j} + G^k_{i-1,j} (\bar U_{real}^{k+1})_{i-1,j} + G^k_{i,j+1} (\bar U_{real}^k)_{i,j+1} + G^k_{i,j-1} (\bar U_{real}^{k+1})_{i,j-1} \nonumber \\
 + \lambda_1(\hat U_{real})_{i,j} +2\lambda_3(\bar U_{real}^k)_{i,j}  \\
	(\mathcal{D}_{real})_{i,j}= \lambda_1 + 2\lambda_3 \left( (\bar U_{real}^k)^2_{i,j} + (\bar U_{im}^k)^2_{i,j} \right)  + G^k_{i+1,j} + G^k_{i-1,j} + G^k_{i,j+1} + G^k_{i,j-1}  \\
 (\mathcal{N}_{im})_{i,j}= \bar G^k_{i+1,j} (\bar U_{im}^k)_{i+1,j} + \bar G^k_{i-1,j} (\bar U_{im}^{k+1})_{i-1,j} + \bar G^k_{i,j+1} (\bar U_{im}^k)_{i,j+1} + \bar G^k_{i,j-1} (\bar U_{im}^{k+1})_{i,j-1} \nonumber \\
 + \lambda_2(\hat U_{im})_{i,j} +2\lambda_3(\bar U_{im}^k)_{i,j} \\
(\mathcal{D}_{im})_{i,j} =	\lambda_2 + 2\lambda_3 \left( (\bar U_{real}^k)^2_{i,j} + (\bar U_{im}^k)^2_{i,j} \right)  + \bar G^k_{i+1,j} + \bar G^k_{i-1,j} + \bar G^k_{i,j+1} + \bar G^k_{i,j-1}
\end{gather}
%
%
%

\section{Experimental results}

In this section, we present  results using both synthetic and real examples. All algorithms presented in this section were run on a computer with an Intel Core i7 2.5 GHz processor, with 16 GB of RAM, and Debian GNU/Linux 9 (stretch) 64-bit operating system installed.  In particular, our model was implemented using the C/ C++ language and the high-performance vector mathematics library Blitz \cite{blitz}. The code and experiments can be found \href{https://github.com/clirlab/wrapped_phase_denoising}{here}.

\subsection{Quality of results}

\begin{figure}[ht]	
	\centering
	\subfigure[]{\includegraphics[width=0.45\textwidth]{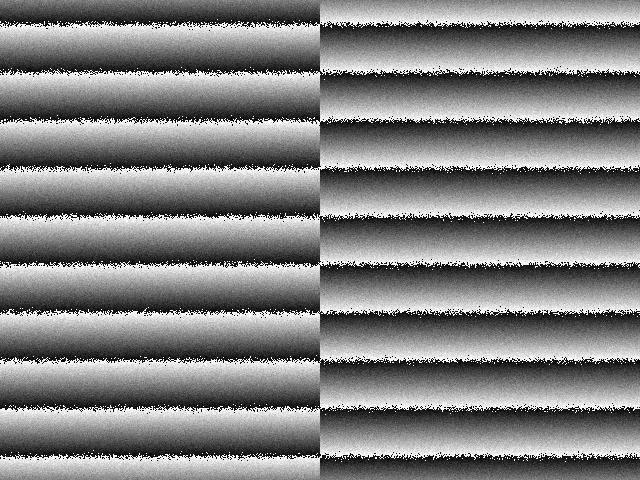}} 
	\subfigure[]{\includegraphics[width=0.45\textwidth]{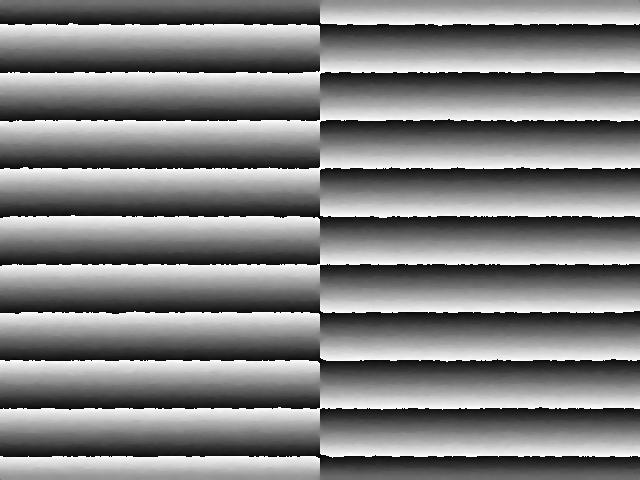}}\\
	\subfigure[]{\includegraphics[width=0.45\textwidth]{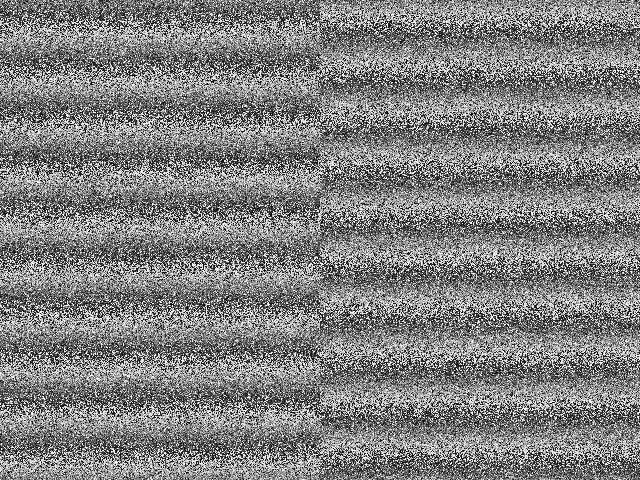}}
	\subfigure[]{\includegraphics[width=0.45\textwidth]{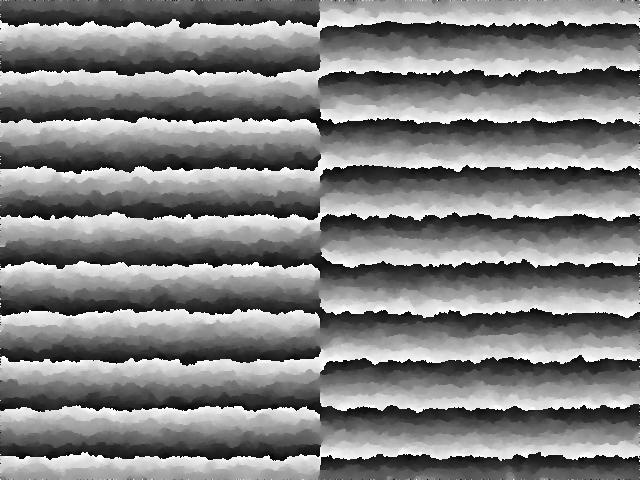}}
	\caption{Left column: noisy wrapped phase map images with 74.41 dB (a) and 43.34 dB (c) SNR levels. Right column: restored phase maps using our model.}
		\label{synt1}
\end{figure}

First, we show results on synthetic images polluted with different levels of additive Gaussian noise. In Figures \ref{synt1}, (a) and (c), we show a wrapped phase image with a vertical discontinuity in the middle with Signal to Noise Ratio (SNR) equal to $74.41$ dB  and $43.34$ dB respectively. Both images were processed by our model resulting in the restored images shown in Figures \ref{synt1}(b) and \ref{synt1}(d). The parameters were fixed to  $\beta=0.001$, $\lambda_1=\lambda_2=2.5$, and $\lambda_3=5$ in both problems. We highlight that noise was fairly removed and the discontinuity preserved.

\begin{figure}[htbp]
	\centering
	\subfigure[]{\includegraphics[width=0.45\textwidth]{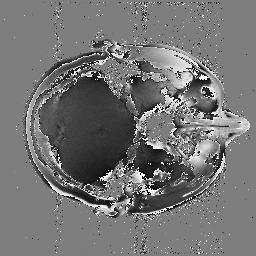}}
	\subfigure[]{\includegraphics[width=0.45\textwidth]{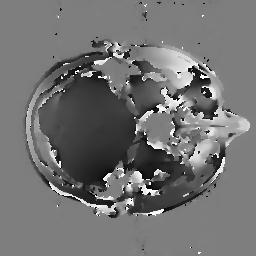}}
	\caption{In (a), a noisy wrapped phase map image from an MRI scan. In (b), the restored wrapped phase map using our model.}
	\label{real1}
\end{figure}

Further, in Figure \ref{real1} we show the noisy and restored image of real problem. Again, restoration is good and edges well preserved. This problem is taken from an MRI scan of the brain.

\subsection{Algorithm performance}
Now, we show the performance of two numerical algorithms for solving our variational model: the gradient descent and the fixed point methods.

In Table \ref{tabla1}, we show the CPU-time and number of iterations taken by each algorithm in solving the problems of Figures \ref{synt1}(a) and (c). For both algorithms, iterations were stopped when the relative residual was less than $\epsilon = 10^{-7}$ with $\epsilon \equiv {||\bar u^{k+1}-\bar u^k||_2}/{||\bar u^k||_2} $

 It can be appreciated from Table \ref{tabla1} that the fixed point is roughly one order of magnitude faster than the gradient descent method as expected.
\begin{table}[h]
	\centering
			\begin{tabular}{lccc}
				Numerical Method & Problem & $\#$ of iterations & CPU-time (s) \\
				\hline
				\hline
				\multirow{2}{*}{Gradient Descent} & Fig. \ref{synt1}(a)& \multicolumn{1}{c}{1,649} & \multicolumn{1}{c}{130} \\\cline{2-4}
				&Fig. \ref{synt1}(c)& \multicolumn{1}{c}{5,159} & \multicolumn{1}{c}{431} \\\cline{2-4} \hline
				\multirow{2}{*}{Fixed Point} &Fig. \ref{synt1}(a)& \multicolumn{1}{c}{104} & \multicolumn{1}{c}{18} \\\cline{2-4}
				&Fig. \ref{synt1}(c)& \multicolumn{1}{c}{307} & \multicolumn{1}{c}{55} \\\cline{2-4}				
              \hline
			\end{tabular}
		\caption{Comparison between gradient descent method and the fixed point method. The later, is much better in terms of CPU-time and number of iterations.}
				\label{tabla1}
	\end{table}
%
%
%
\subsection{Comparison against other denoising algorithms}
A natural question would be to know whether state of the art algorithms for image denoising deliver a good restoration of $\bar U_{real}$ and $\bar U_{im}$ and intrinsically the wrapped phase map. To this end, we tested BM3D \cite{dabov2007bm3d}, Non Local Means (NLM) \cite{buades2005nlm},  and also the de facto method for denoising wrapped phase maps Str\"obel's method \cite{strobel1996processing} on the synthetic problems of Figures \ref{synt1}(a) and (c) and the real problem of Figure \ref{real1}(a). 

To compare results, we use two metrics: the Mean Squared Error (MSE) and the Universal Image Quality Index (IQI) defined below for any variable $y \in \mathbb{R}^{m \times n}$. 
\begin{gather}
MSE = \frac{1}{m*n}{\sum_{i=1}^m\sum_{j=1}^n(y(i,j)-\bar y(i,j))^2},
\end{gather}
and
\begin{gather}
	IQI = 1-\frac{\sum_{i=1}^m\sum_{j=1}^n \left (y(i,j)-\bar y(i,j) \right)^2}{\sum_{i=1}^m\sum_{j=1}^n y(i,j)^2}.
\end{gather}
The dynamic range of the IQI is $[-1, 1]$, where the worst value is $-1$ and the best value is $1$.

\begin{figure}[htbp]
	\centering
	\subfigure[]{\includegraphics[width=0.45\textwidth]{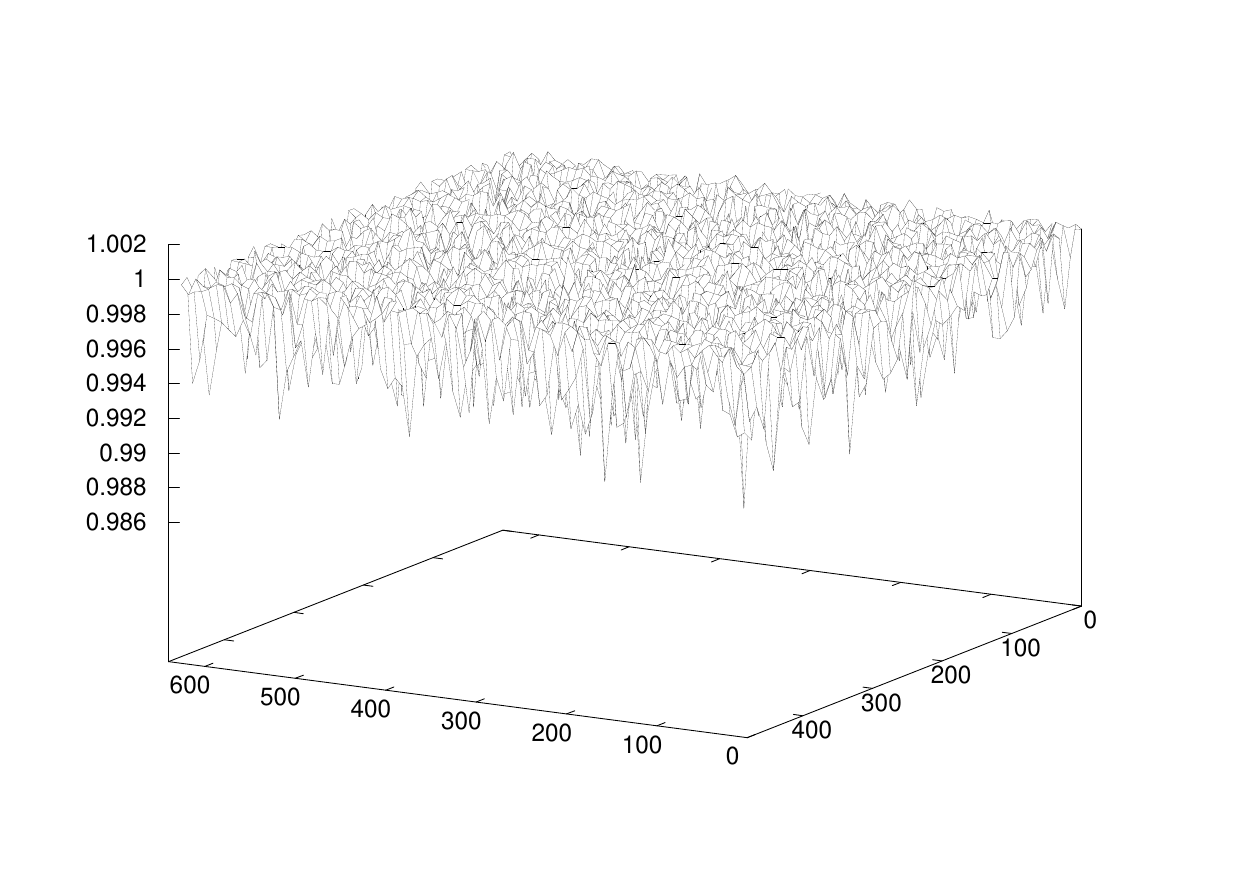}}
	\subfigure[]{\includegraphics[width=0.45\textwidth]{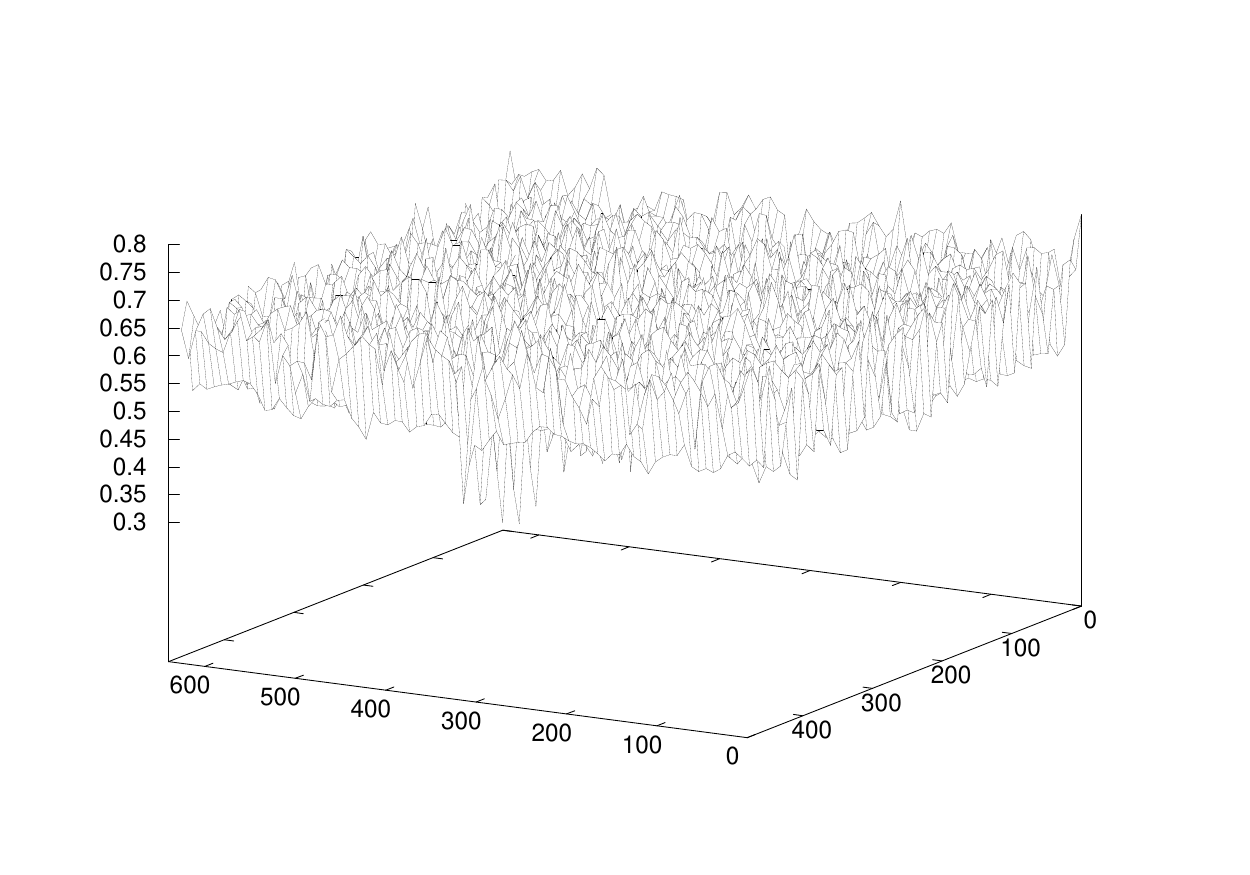}}
	\subfigure[]{\includegraphics[width=0.45\textwidth]{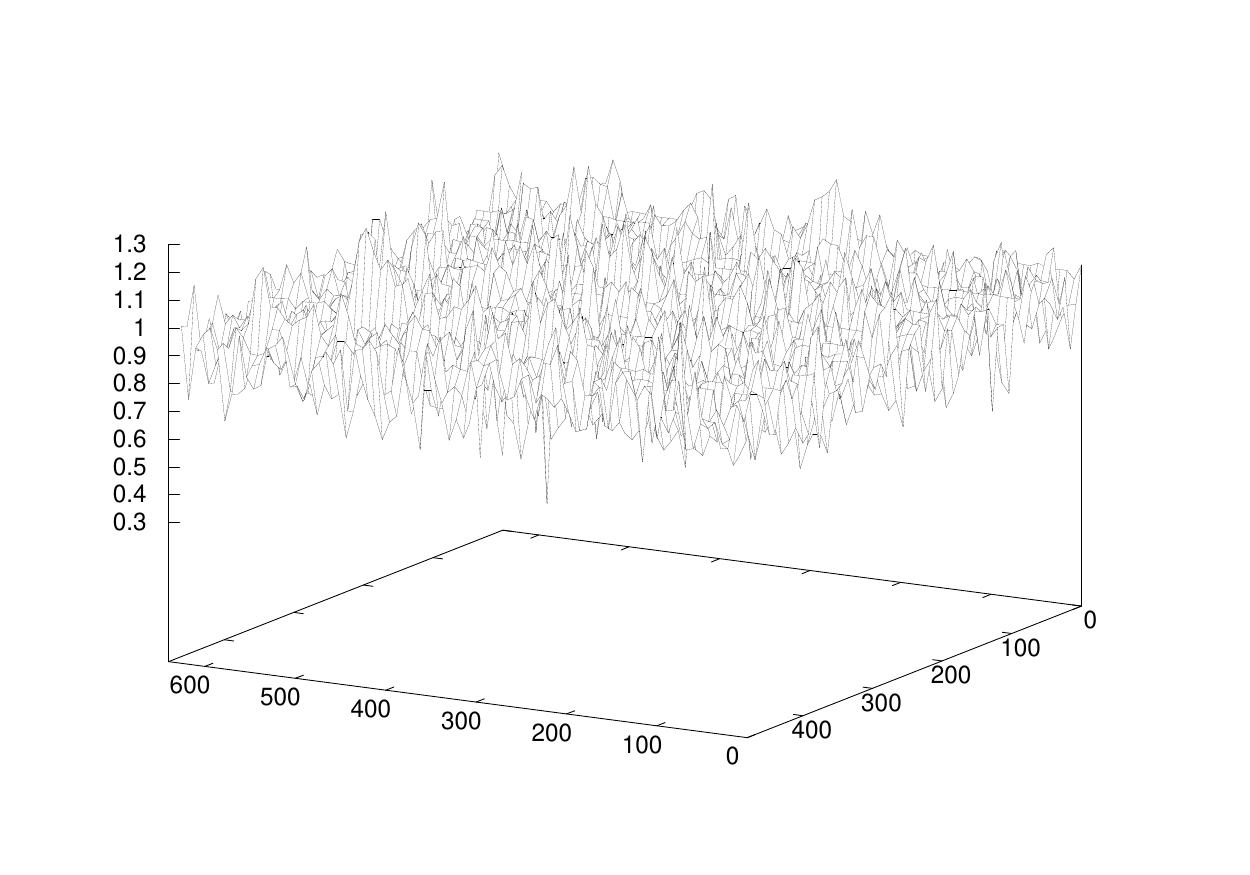}}
	\subfigure[]{\includegraphics[width=0.45\textwidth]{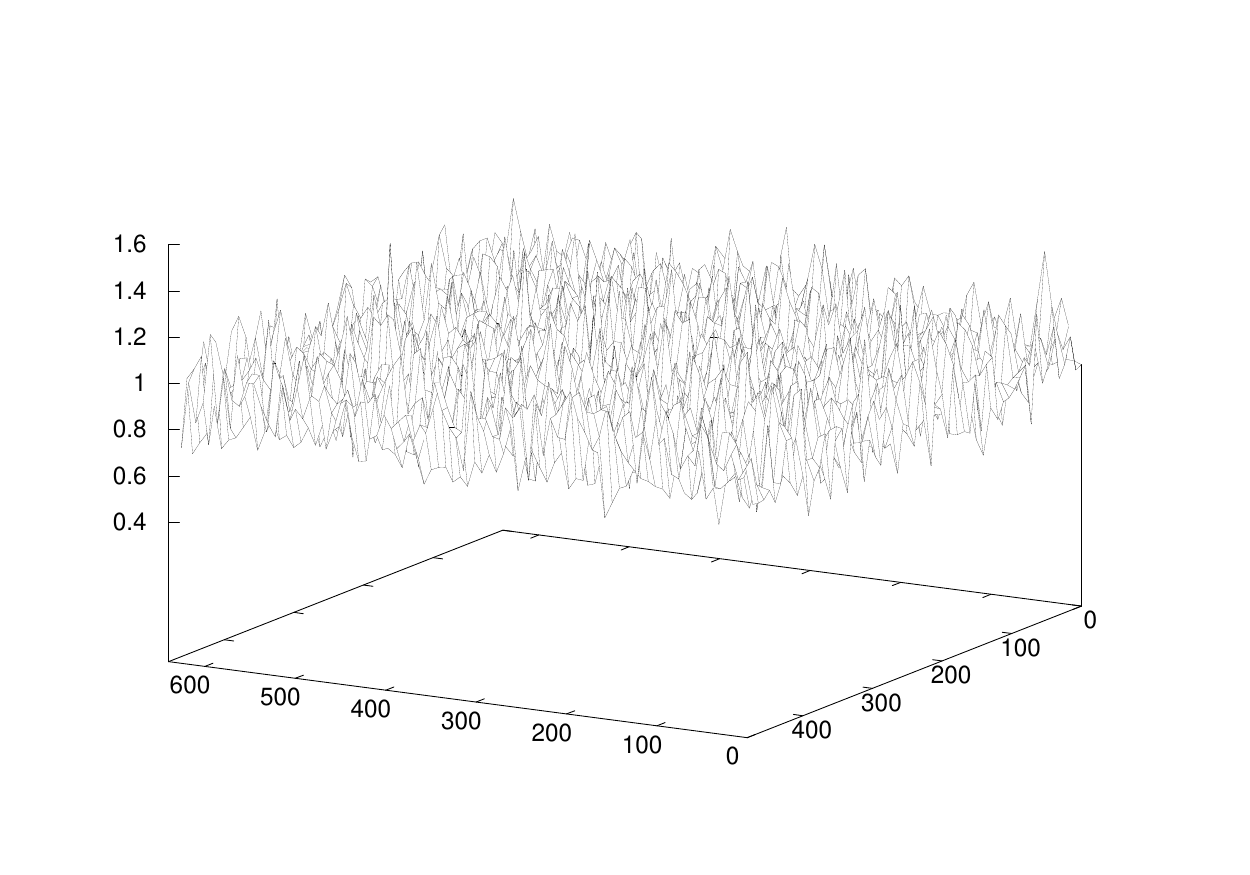}}
	\caption{Plot of $||\bar U_{real}^2 + \bar U_{im}^2 ||_2$ for the problem in Figure \ref{synt1}(a). Notice in (a) that our model is able to preserve the trigonometric relation between $\bar U_{real}^2 + \bar U_{im}^2 =1$. On the other hand, Str\"obel's (b) method, BM3D (c) and NLM (d) do not preserve the Pythagorean property.}
	\label{fig-comp-PF-Strobel-bm3d-nlmeans}
\end{figure}

\begin{figure}[htbp]
	\centering
	\subfigure[]{\includegraphics[width=0.22\textwidth]{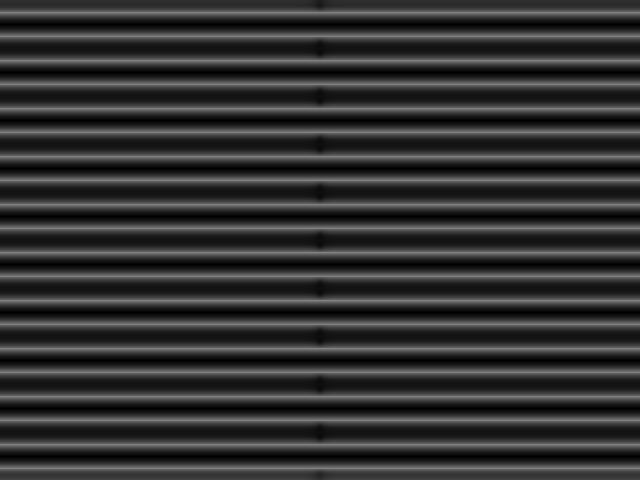}}
	\subfigure[]{\includegraphics[width=0.22\textwidth]{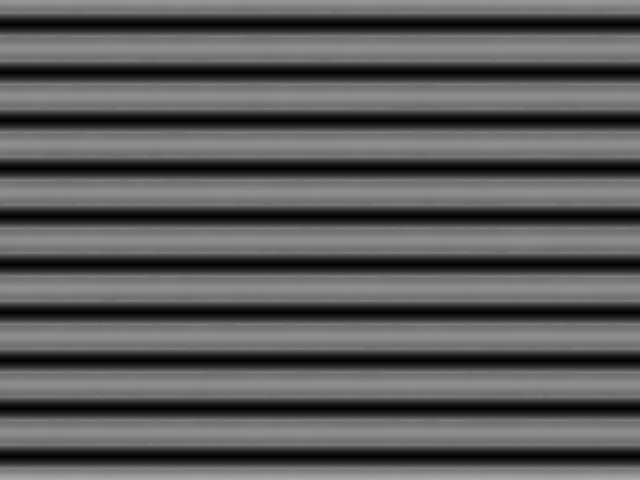}}
	\subfigure[]{\includegraphics[width=0.22\textwidth]{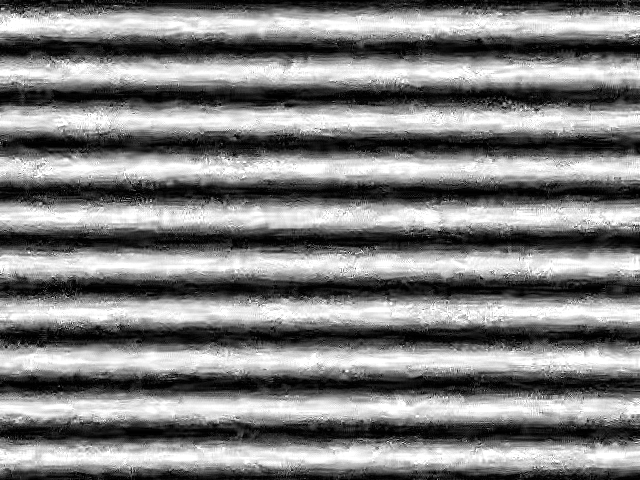}}
	\subfigure[]{\includegraphics[width=0.22\textwidth]{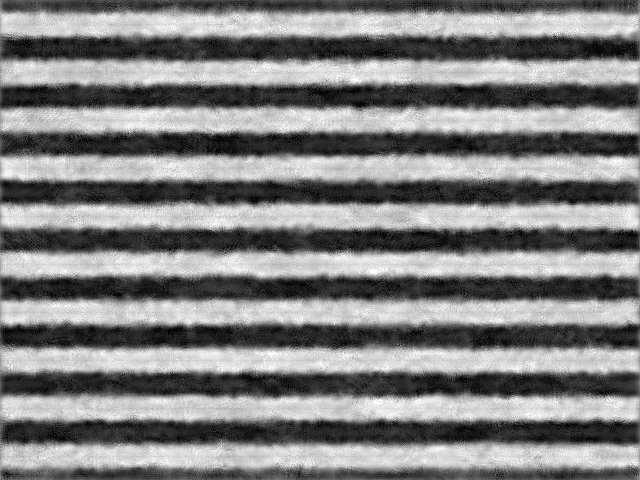}}
	\subfigure[]{\includegraphics[width=0.22\textwidth]{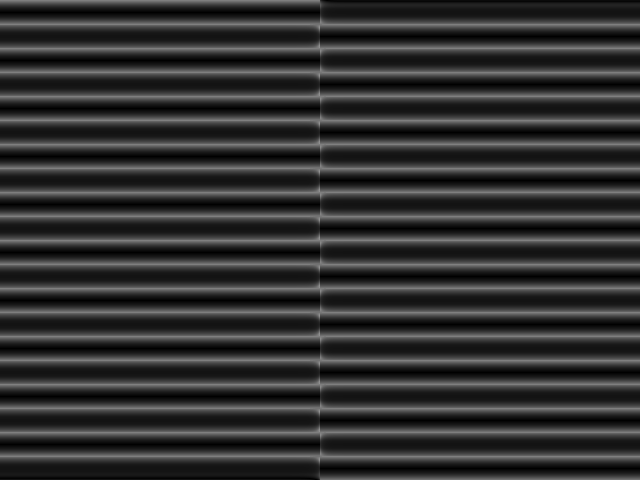}}
	\subfigure[]{\includegraphics[width=0.22\textwidth]{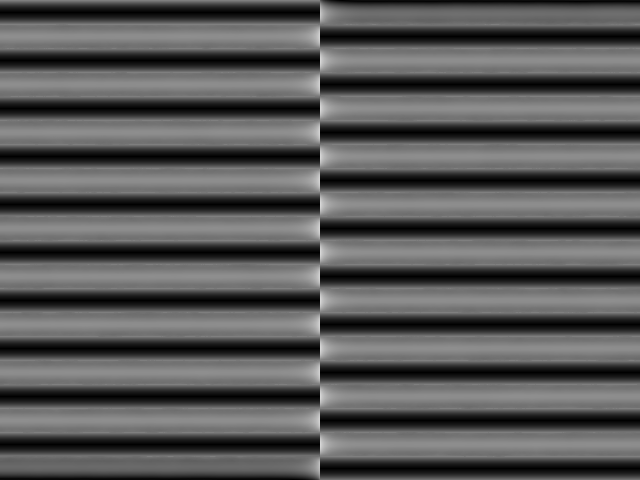}}
	\subfigure[]{\includegraphics[width=0.22\textwidth]{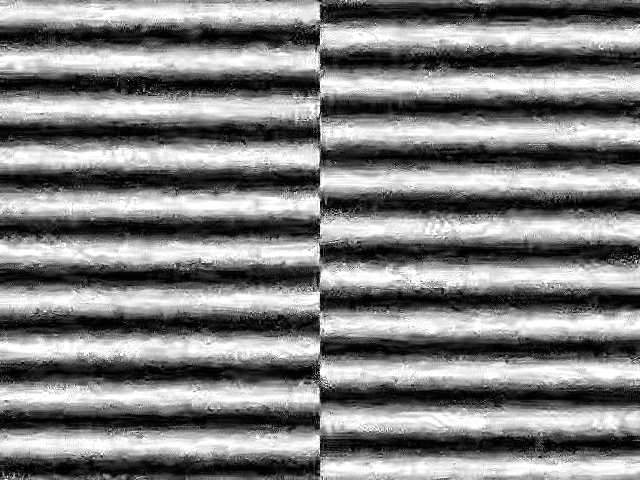}}
	\subfigure[]{\includegraphics[width=0.22\textwidth]{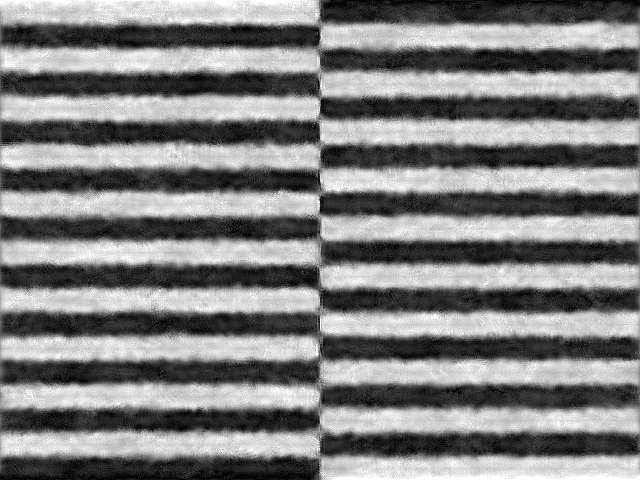}}
	\caption{In the top row, we show the difference $|\bar{U}_{real}-\hat{U}_{real}|$ for the problem in Figure \ref{synt1}(a) for all models tested in this work. The bottom row shows the difference $|\bar{U}_{im}-\hat{U}_{im}|$ for the same models. Notice that in our model, first column, there is only a small error that shows as slim horizontal lines. On the other hand, Str\"obel's result in the second column, BM3D, and NLM results in the third and fourth columns respectively, show a larger error present at many other points in the domain.}
	\label{fig-comp-Ic-sintetica-PF-Strobel-bm3d-nlmeans}
\end{figure}

The obtained results from these two metrics are presented in Table \ref{tab-metricas}. By far, our model gets a much better restoration. Even Strobel's method overcomes BM3D and NLM illustrating that for denoising a wrapped phase map image, the relation between $\bar U_{real}$ and $\bar U_{im}$ has to be preserved. This fact is more evident if we look at Figures \ref{fig-comp-PF-Strobel-bm3d-nlmeans} and \ref{fig-comp-Ic-sintetica-PF-Strobel-bm3d-nlmeans} showing that enforcing the Pythagorean trigonometric identity between the real and imaginary parts of the phase map effectively improves the wrapped phase map restoration reducing this way inaccuracies in possible applications.
\begin{table}[ht!]
	\centering
	\caption{Comparison among different denoising algorithms.}
	\begin{tabular}{c|cc|cc}
		Model &  \multicolumn{2}{c}{MSE} &  \multicolumn{2}{c}{IQI} \\ \hline  \hline \\
		& $\bar{U}_{real}$ & $\bar{U}_{im}$ & $\bar{U}_{real}$ & $\bar{U}_{im}$\\ \hline
		Our model & 0.0555 & 0.0542 & 0.8515 & 0.8551 \\
		Str\"obel & 0.1537 & 0.1489 & 0.5891 & 0.6021\\
	BM3D & 0.3696 & 0.3687 & 0.0125 & 0.0149  \\
		NLM & 0.3478 & 0.3511 & 0.0706 & 0.0621 \\
	\end{tabular}
	\label{tab-metricas}
\end{table}
\section{Conclusions}
In this paper, we introduced a new variational model for denoising wrapped phase maps. The two main properties of our model are: regularization is done through total variation therefore the model preserves discontinuities and the fundamental Pythagorean trigonometric identity between the real and imaginary parts of the phase map is enforced improving considerably the restoration of the phase map. Further, the existence and uniqueness of the solution of our model is analyzed and a fast fixed point algorithm for the numerical realization is provided. Theoretical analysis of this fixed point method to prove its convergence is also provided. Finally, restoration results on synthetic and real examples are given as evidence of the quality of restoration and comparisons against state of the art methods also discussed.



\section{Declarations}
The authors of this manuscript declare that there are no conflicts of interest to report. All co-authors participated in the writing of the manuscript and agreed with its contents. We certify that the work submitted is original and is not under review at any other publication.
\bibliographystyle{apalike}

\end{document}